\documentclass[a4paper,11pt]{article}

\title{{\bf Connectivity Threshold of Random Geometric Graphs with Cantor Distributed Vertices}}
\author{Antar Bandyopadhyay\footnote{E-Mail: {\tt antar@isid.ac.in}} \\  Farkhondeh Sajadi\footnote{E-Mail: {\tt sajadi7r@isid.ac.in}} \\  
\ \\
\ \\
Theoretical Statistics and Mathematics Unit\\
Indian Statistical Institute, Delhi Centre\\
7 S. J. S. Sansanwal Marg \\
New Delhi 110016 \\ 
INDIA}
\date{April 2, 2012}

\usepackage{amsfonts,amsmath,amssymb,mathrsfs,amsthm}
\usepackage{bm}
\usepackage{graphicx,graphics,subfig}
\theoremstyle{plain}
\newtheorem{Theorem}{Theorem}

\newcommand{\bP}{\mathbb{P}}

\begin{document}

\maketitle

\begin{abstract}
For connectivity of \emph{random geometric graphs}, where there is no density for underlying distribution of the vertices, we consider 
$n$ i.i.d. \emph{Cantor} distributed points on $[0,1]$.  We show that for this random geometric graph, the connectivity threshold $R_{n}$, 
converges almost surely to a constant $1-2\phi$ where $0 < \phi < 1/2$, 
which for the standard Cantor distribution  is $1/3$. We also show that 
$\left\| R_n - \left(1 - 2 \phi \right) \right\|_1 \sim 2 \, C\left(\phi\right) \, n^{-1/d_{\phi}}$
where $C\left(\phi\right) > 0$ is a constant and $d_{\phi} := - {\log 2}/{\log \phi}$ is the
\emph{Hausdorff dimension} of the generalized Cantor set with parameter $\phi$. \\


\end{abstract}

\section{Introduction}
\label{Sec:Intro}
\subsection{Background and motivation}
\label{SubSec:Back and Moti}
A \emph{random geometric graph} consists of a set of vertices, distributed randomly over some metric space,  
in which two distinct such vertices are joined by an edge, if the distance between them is sufficiently small. More precisely, 
let $V_{n}$ be a set of $n$ points in $\mathbb{R}^d$, distributed independently according to some distribution $F$ on $\mathbb{R}^d$. 
Let $r$ be a fixed positive real number. Then, the random geometric graph $\mathcal{G}=\mathcal{G}(V_{n}, r)$ is a graph with vertex set $V_{n}$ 
where two vertices $\mathbf{v} = (v_{1}, \ldots , v_{d})$ and $\mathbf{u} = (u_{1}, \ldots , u_{d})$ in $V_{n}$ are adjacent if and only 
if $\left\|\mathbf{v}-\mathbf{u}\right\| \leq r$ where $\left\|.\right\|$ is some norm on $\mathbb{R}^d$.

A considerable amount of work has been done on the \emph{connectivity threshold} defined as
\begin{equation}
R_{n}=\inf \left\{ r>0  \,\Big\vert\,  \mathcal{G}(V_{n},r)~ \text{is connected} \,\right\} \,.
\label{Defn of R}
\end{equation}

The case when the vertices are assumed to be uniformly distributed in $[0,1]^d$, \cite{ApRu02} showed that with probability one
\begin{equation}
\lim_{n\rightarrow\infty} \frac{n}{\log n} \, R^{d}_{n} = \left\{
\begin{array}{rl}
1 & \text{for } d=1,\\
\frac{1}{2d} & \text{for } d\geq 2 \quad 
\end{array} \right. \,,
\label{Eq:Unif-Almost-Sure-Limit}
\end{equation}
when the norm $\left\|.\right\|$ is taken to be the ${\mathcal L}_{\infty}$ or the $\sup$ norm. Later \cite{Pen03} 
showed that the limit in
\eqref{Eq:Unif-Almost-Sure-Limit} holds but with different constants for any ${\mathcal L}_p$ norm for $1 \leq p \leq \infty$.  
\cite{Pen99} considered the case when the distribution $F$ has a continuous density $f$ with respect to the Lebesgue measure
which remains bounded away from $0$ on the support of $F$. Under certain technical assumptions such as smooth boundary for the
support, he showed that with probability one, 
\[
\lim_{n \rightarrow \infty} \frac{n}{\log n} \, R^{d}_{n} = C
\]
where $C$ is an explicit 
constant which depends on the dimension $d$ and essential infimum of $f$ and its value on the boundary of the support. Recently, 
\cite{SaSa10} [personal communication], studied a case when the density $f$ 
of the underlying distribution may have minimum zero. They in particular, proved that when the support of $f$ is $[0, 1]$ and $f$ is bounded below on 
any compact subset not containing the origin but it is regularly varying at the origin, then $R_n/F^{-1}(1/n)$ has a weak limit. 

The proof by \cite{SaSa10} can
easily be generalized to the case where the density is zero at finitely many points. A question then naturally arises what happens to the case when the 
distribution function is flat on some intervals, that is, if density exists then it will be zero on some intervals. Also what happens in the
somewhat extreme case, when the density may not exist even though the distribution function is continuous and has flat parts. To consider these
questions, 
in this paper we study the connectivity of random geometric graphs where the underlying distribution of the 
vertices has no mass and is also singular with respect to the Lebesgue measure, that is, it has no density. 
For that, we consider the 
\emph{generalized Cantor distribution} with parameter $\phi$ denoted by $Cantor(\phi)$ 
as the underlying distribution of the vertices of the graph. The distribution function is then flat on 
infinitely many intervals. 
We will show that the connectivity threshold converges almost surely to the length of the longest
flat part of the distribution function and we also provide some finer asymptotic of the same.

\subsection{Preliminaries}
\label{SubSec:Prelim}
In this subsection, we discuss the \emph{Cantor set} and \emph{Cantor distribution} which is defined on it. 

\subsubsection{Cantor Set}
The Cantor set was first discovered by  \cite{Sm1875} but became popular after \cite{Can1883}. 
The Standard Cantor set is constructed on the interval $[0,1]$ as follows. One successively removes the open middle third of each subinterval of the 
previous set. 
More precisely, starting with $C_0 := [0,1]$, we inductively define
\[
C_{n+1} := \bigcup_{k=1}^{2^n} 
            \left( \left[a_{n,k}, a_{n,k} + \frac{b_{n,k} - a_{n,k}}{3} \right] 
                   \cup 
                   \left[b_{n,k} - \frac{b_{n,k} - a_{n,k}}{3}, b_{n,k} \right] \right)
\]
where  
$C_n := \displaystyle{\mathop{\cup}\limits_{k=1}^{2^n} \left[a_{n,k} , b_{n,k}\right]}$.
The Standard Cantor set is then defined as 
$C=\bigcap_{n=o}^{\infty} C_n$.
It is known 
that the Hausdorff dimension of the standard Cantor set is $\frac{\log 2}{\log 3}$ (see Theorem 2.1 of Chapter 7 of \cite{StSh05}).

For constructing the generalized Cantor set, we start with unit interval $[0,1]$ and at 
first stage, we delete the interval $(\phi, 1-\phi)$ where $0 < \phi < 1/2$. Then, this procedure is reiterated with the two segments 
$[0,\phi]$ and $[1-\phi,1]$. We continue ad infinitum. The Hausdorff dimension of this set is given by
$d_{\phi} := - \frac{\log 2}{\log \phi}$ (see Exercise 8 of Chapter 7 of \cite{StSh05}). 
Note that the standard Cantor set is a special case when $\phi = 1/3$.

\subsubsection{Cantor distribution}
The \emph{Cantor distribution} with parameter $\phi$ where $0 < \phi < 1/2$ is the distribution  
of a random variable $X$ defined by 
\begin{equation} X=\sum_{i=1}^{\infty}\phi^{i-1} Z_{i} \label{Defn Cantor}\end{equation}
where $Z_{i}$ are i.i.d. with $\mathbb{P}[Z_{i}=0]= \mathbb{P}[Z_{i}=1-\phi]=1/2$. 
If a random variable $X$ admits a representation of the form \eqref{Defn Cantor} then we will say that $X$ has a 
Cantor distribution with parameter $\phi$, and write $X \sim \text{Cantor}(\phi)$. 
Observe that $\text{Cantor}\left(\phi\right)$ is self-similar, in the sense that, 
\begin{equation}
X \stackrel{d}= \left\{
\begin{array}{lr}
\phi X & \text{with probability}~ 1/2\\
\phi X + 1 - \phi & \text{with probability}~ 1/2
\end{array} \right.
\label{Equ:Fundamental-Recursion}
\end{equation}
This follows easily by conditioning on $Z_{1}$. It is worth noting here that if $X \sim \text{Cantor}(\phi)$ then
so does $1 - X$. 

Note that for $\phi=1/3$ we obtain the \emph{standard Cantor distribution}.

\section{Main results}
\label{Sec:Main Result}
Let $X_{1}, X_{2}, ..., X_{n}$ be independent and identically distributed random variables with $\text{Cantor}(\phi)$ distribution on $[0,1]$.  
Given the graph $\mathcal{G}=\mathcal{G}(V_{n}, r)$, where $V_{n}=\{X_{1}, X_{2}, \ldots, X_{n}\}$, let $R_{n}$ be defined as in \eqref{Defn of R}.
\begin{Theorem}
\label{Thm:Convergence}
For any $0 < \phi < 1/2$, as $n \longrightarrow \infty$ we have
\begin{equation}
R_{n} \longrightarrow 1-2 \phi \mbox{\ a.s.}
\end{equation}
\end{Theorem}

Our next theorem gives finer asymptotics but before we state the theorem, we provide here some basic notations and facts. 
Let $m_{n}:=\!\min\{X_{1}, X_{2}, \ldots, X_{n}\!\}$. 
Using \eqref{Equ:Fundamental-Recursion} we get
\begin{equation}
m_{n} \stackrel{d}= \left\{
\begin{array}{lcl}
\phi m_{k} & \text{with probability} &  2^{-n}\binom{n}{k} \text{ for } k=1,2,...,n\\
\phi m_{n}+1-\phi & \text{with probability} & 2^{-n}   \label{Rec for Min order}
\end{array} \right.
\end{equation}
Let  $a_{n}:=\mathbb{E}[m_{n}]$. Using \eqref{Rec for Min order} \cite{Hos94} derived the following recursion 
formula for the sequence $\left(a_n\right)$
\begin{equation}
\left(2^n-2\phi\right) a_{n}=1-\phi+\phi \sum_{k=1}^{n-1} \binom{n}{k} a_{k}, \quad n\geq1    \label{Rec for an} 
\end{equation}
Moreover \cite{KnPr96} showed that whenever $0 < \phi < 1/2$ then as $n \rightarrow \infty$,
\begin{equation}
\frac{a_n}{n^{-\frac{1}{d_{\phi}}}} \longrightarrow C\left(\phi\right)   \,,
\label{Asym for an}
\end{equation}
where 
\begin{equation}
C\left(\phi\right) := \frac{(1-\phi)(1-2\phi)}{\phi \log 2} \Gamma(-\log_{2} \phi)\zeta(-\log_{2} \phi) \,,
\label{Equ:An-asym-Cons}
\end{equation}
and $d_{\phi} = - \frac{\log 2}{\log \phi}$ is the Hausdorff dimension of the generalized Cantor set. 
Here $\Gamma(\cdot)$ and $\zeta(\cdot)$ are the Gamma and Riemann zeta functions, respectively. 

Our next theorem gives the rate convergence of $R_n$ to $\left(1 - 2 \phi\right)$ in terms of the ${\mathcal L}_1$ norm.  
\begin{Theorem}
\label{Thm:Rate-of-R}
For any $0 < \phi < 1/2$, 
as $n \longrightarrow \infty$ we have 
\begin{equation}
\frac{\left\|R_{n}-(1-2\phi) \right\|_{1}}{n^{- \frac{1}{d_{\phi}}}} \longrightarrow 2 C\left(\phi\right)\,,
\label{Equ:Rate}
\end{equation}
where $C\left(\phi\right)$ is as in equation \eqref{Equ:An-asym-Cons} and 
$\left\| \cdot \right\|_1$ is the ${\mathcal L}_1$ norm. 
\end{Theorem}

\section{Proof of the theorems}
\label{Sec:Proofs}
\subsection{Proof of Theorem \ref{Thm:Convergence}}
\label{SubSec:The1}
We draw a sample of size $n$ from $Cantor(\phi)$ on $[0,1]$. Let $N_n$ be the number of elements falling in the 
subinterval $[0,\phi]$ and $n-N_n$ in $[1-\phi,1]$. From the construction $N_{n} \sim \text{Bin}\left(n,\frac{1}{2}\right)$.
In selecting this sample of size $n$, there are three cases which may happen. Some of these points may fall in interval $[0,\phi]$ and rest in interval 
$[1-\phi,1]$. That means $N_{n} \notin \{0,n\}$. In this case the distance between the points in $[0,\phi]$ and $[1-\phi,1]$ is at least $1 - 2\phi$. 
The other cases are when all points fall in $[0,\phi]$ or all fall in $[1-\phi,1]$ , which in this case $N_{n}=n$ or $N_{n}=0$.
Let $\displaystyle{m_{n}=\min_{1\leq i \leq n} X_{i}}$, $\displaystyle{M_{n}=\max_{1\leq i \leq n} X_{i}}$ and we define
\begin{equation}
L_{n}:= \max \left\{X_{i} |~ 1\leq i\leq n ~ \text{and}~ X_{i} \in [0,\phi]  \right\}
\end{equation}
and
\begin{equation}
U_{n}:= \min \left\{X_{i} |~ 1\leq i\leq n ~ \text{and}~ X_{i} \in [1-\phi,1] \right\} \,.
\end{equation}
We will take $L_n = 0$ (and similarly $U_n = 0$) if the corresponding set is empty. 

Now find a $K \equiv K\left(\phi\right)$ such that $\phi^K < \frac{1}{2}\left(1 - \phi\right)\left(1 - 2 \phi\right)$. 
Note that such a $K < \infty$ exists since $0 < \phi < 1$. 
Let $I_1, I_2, \ldots, I_{2^K}$ be the $2^K$ sub-intervals of length $\phi^K$ which are part of the $K^{\mbox{th}}$ stage of the
``removal of middle interval'' for obtaining the generalized Cantor set with parameter $\phi$. For $1 \leq j \leq 2^K$ define
$N_j := \sum_{i=1}^n \bm{1}\left(X_i \in I_j\right)$, which is the number of sample points in the sub-interval $I_j$. From the construction
of the the generalized Cantor distribution with parameter $\phi$ it follows that
\begin{equation}
\bm{N}_K := \left(N_1, N_2, \ldots, N_{2^K}\right) \sim \mbox{Multinomial}\left(n; \left(\frac{1}{2^K}, \frac{1}{2^K}, \cdots, \frac{1}{2^K}\right) \right)\,,
\label{Equ:Multinimial}
\end{equation}
and $N_n = \displaystyle{\mathop{\sum}\limits_{I_j \subseteq \left[0,\phi\right]} N_j}$. 
Consider the event $E_n := \displaystyle{\mathop{\cap}\limits_{j=1}^{2^k} \left[ N_j \geq 1 \right]}$. 
Observe that on the event $E_n$ the maximum inter point distance between two points in $[0,\phi]$ as well as in $[1-\phi, 1]$ is at most
$2 \phi^K + \phi \left(1 - 2\phi\right) < 1 - 2 \phi$ by the choice of $K$. Thus on $E_n$
we must have $R_n = U_n - L_n$ and so we can write
\begin{equation} 
R_{n}=\left(U_{n}-L_{n}\right) \, \bm{1}_{E_n} + R_n^* \, \bm{1}_{E_n^c}
\label{Rn} 
\end{equation}
where $R_n^*$ is a random variable such that $0 < R_n^* < \phi$ a.s. 

Observe that conditioned on $\left[N_1 = r_1, N_2 = r_2, \ldots, N_{2^k}= r_{2^k} \right]$ we have 
$U_n \stackrel{d}= 1 - \phi + \phi m_{n-k}$ and $L_n \stackrel{d}= \phi M_{k}$ and $N_n = k$
where $k = \displaystyle{\mathop\sum\limits_{I_j \subseteq \left[0,\phi\right]} r_j}$. 
More generally
\begin{equation}
\left(\left(L_n, U_n\right), \bm{N}_K \right)_{n \geq 1} 
\stackrel{d}= 
\left( \left(  \phi M_{N_n}, 1 - \phi + \phi m_{n-N_n} \right), \bm{N}_K \right)_{n \geq 1} \,.
\label{Equ:process-equality} 
\end{equation}
Note that to be technically correct we define $M_0 = m_0 = 0$. 

Now it is easy to see that $m_n \longrightarrow 0$ and $M_n \longrightarrow 1$ a.s. But
by the SLLN,  $N_n/n \longrightarrow 1/2$ a.s., thus
both $\left(N_n\right)$ and $\left(n - N_n\right)$ are two subsequences which are converging to infinity a.s.
Moreover
\begin{equation}
\bP\left(E_n^c\right) \leq \sum_{j=1}^{2^K} \bP\left(N_j = 0 \right) 
                      =  2^K \left(1 - \frac{1}{2^K}\right)^n 
                      =  2^K \exp\left(- \alpha_K n \right) 
\label{Equ:Exponential-Decay} \,,
\end{equation}
where $\alpha_K = - \log \left(1 - \frac{1}{2^K}\right) > 0$. Thus $\sum_{n=1}^{\infty} \bP\left(E_n^c\right) < \infty$, so by the
First Borel-Cantelli Lemma we have 
\[
\bP\left(E_n^c \mbox{\ infinitely often\ }\right) = 0 \,\, \Rightarrow \,\, \bP\left(E_n \mbox{\ eventually\ }\right) = 1 \,.
\] 
In other words $\bm{1}_{E_n} \longrightarrow 1$ a.s. and $\bm{1}_{E_n^c} \longrightarrow 0$ a.s. 
Finally observing that $0 \leq R_n^* \leq \phi$ we get from equations \eqref{Rn} and \eqref{Equ:process-equality} 
\[ 
R_n \longrightarrow \left(1 - 2 \phi \right) \,.
\]
\ \ \hfill $\square$

\subsection{Proof of Theorem \ref{Thm:Rate-of-R}}
\label{SubSec:The2}
We start by observing
\begin{eqnarray}
&   & \mathbb{E}\left[\left\vert R_{n} - \left(1 - 2 \phi\right) \right\vert \right] \nonumber \\
& = & \mathbb{E}\left[\left( R_{n} - \left(1 - 2 \phi\right) \right) \bm{1}_{E_n}\right] + 
      \mathbb{E}\left[\left\vert R_{n}^* - \left(1 - 2 \phi\right) \right\vert \bm{1}_{E_n^c} \right] \nonumber \\
& = & \mathbb{E}\left[\left( U_n - L_n - \left(1 - 2 \phi\right) \right) \bm{1}_{E_n} \bm{1}_{K \leq N_n \leq n-K} \right] + 
      \mathbb{E}\left[\left\vert R_{n}^* - \left(1 - 2 \phi\right) \right\vert \bm{1}_{E_n^c} \right] \nonumber \\
& = & \mathbb{E}\left[\left( U_n - L_n - \left(1 - 2 \phi\right) \right) \bm{1}_{K \leq N_n \leq n-K} \right]  \nonumber \\ 
&   & - 
      \mathbb{E}\left[\left( U_n - L_n - \left(1 - 2 \phi\right) \right) \bm{1}_{E_n^c} \bm{1}_{K \leq N_n \leq n-K} \right] + \mathbb{E}\left[\left\vert R_{n}^* - \left(1 - 2 \phi\right) \right\vert \bm{1}_{E_n^c} \right] \nonumber \\
& = & \mathbb{E}\left[\left( U_n - L_n - \left(1 - 2 \phi\right) \right) \bm{1}_{1 \leq N_n \leq n-1} \right]  \nonumber \\ 
&   & - 
      \mathbb{E}\left[\left( U_n - L_n - \left(1 - 2 \phi\right) \right) \bm{1}_{E_n^c} \bm{1}_{1 \leq N_n \leq n-1} \right] + \mathbb{E}\left[\left\vert R_{n}^* - \left(1 - 2 \phi\right) \right\vert \bm{1}_{E_n^c} \right] \,. \nonumber\\ \label{Equ:L1-Norm}
\end{eqnarray}
In above the first equality holds because of \eqref{Rn} and the fact that on the event $E_n$ we must have $R_n > 1 - 2 \phi$. Second, third
and forth equalities follows from the simple fact that $E_n \subseteq \left[ K \leq N_n \leq n-K \right]$. 

Now recall that $a_{n}=\mathbb{E}[m_n]$ so for the first part of the right-hand side of the equation \eqref{Equ:L1-Norm} we can write
\begin{eqnarray}
&   & \mathbb{E}\left[\left( U_n - L_n - \left(1 - 2 \phi\right) \right) \bm{1}_{1 \leq N_n \leq n-1} \right] \nonumber \\
& = & \frac{\phi}{2^n} \sum_{k=1}^{n-1} \binom{n}{k} (a_{n-k}+a_{k}) \nonumber \\
& = & \frac{1}{2^{n-1}}\left[ \left(2^n - 2 \phi\right) a_n - \left(1-\phi\right) \right] \,, \label{Equ:L1-Norm-I}
\end{eqnarray}
where the last equality follows from \eqref{Rec for an}.
The other two parts of the right-hand side of the equation \eqref{Equ:L1-Norm} 
are bounded in absolute value by 
\[ 
\bP\left(E_n^c\right) \leq 2^K \exp\left(- \alpha_K n \right) 
\]
because of \eqref{Equ:Exponential-Decay}. 
Now observe that from equation \eqref{Asym for an} we get that $a_n \sim C\left(\phi\right) n^{-\frac{1}{d_{\phi}}}$ where
$d_{\phi} = - \frac{\log 2}{\log \phi}$ is the Hausdorff dimension of the generalized Cantor set. Thus using \eqref{Equ:L1-Norm} and 
\eqref{Equ:L1-Norm-I} we conclude that
\[
\frac{\mathbb{E}\left[ \left\vert R_{n}-(1-2\phi) \right\vert \right]}{a_n} \longrightarrow 2  \quad \text{as} \quad n\longrightarrow \infty \,.
\]
This completes the proof using \eqref{Asym for an}.  
\ \ \hfill $\square$

\section{Final Remarks}
It is worth noting here that our proofs depend on the recursive nature of the generalized Cantor ditsribution 
(see equation \eqref{Equ:Fundamental-Recursion}). Thus unfortunately, they 
do not have obvious extensions for other singular distributions. It will be interesting to derive a version of Theorem \ref{Thm:Convergence}
for a general singular distribution with no mass and flat parts. Intuitively it seems that the final limit should be the length of the longest flat part.
It will be more interesting if Theorem  \ref{Thm:Rate-of-R} can also be generalized for general singular distributions
with no mass and flat parts where $(1-2\phi)$ is replaced by the
length of the longest flat part and $d_{\phi}$ is replaced by the Hausdorff dimension of the support.


\bibliographystyle{unsrt}
\bibliography{reference}

\begin{thebibliography}{1}

\bibitem{ApRu02}
Martin J.~B. Appel and Ralph~P. Russo.
\newblock The connectivity of a graph on uniform points on {$[0,1]^d$}.
\newblock {\em Statist. Probab. Lett.}, 60(4):351--357, 2002.

\bibitem{Pen03}
Mathew Penrose.
\newblock {\em Random geometric graphs}, volume~5 of {\em Oxford Studies in
  Probability}.
\newblock Oxford University Press, Oxford, 2003.

\bibitem{Pen99}
Mathew~D. Penrose.
\newblock A strong law for the largest nearest-neighbour link between random
  points.
\newblock {\em J. London Math. Soc. (2)}, 60(3):951--960, 1999.

\bibitem{SaSa10}
A.~Sarkar and B.~Saurabh.
\newblock Connectivity distance in random geometric graph.
\newblock {\em Unpublished results}, 2010.

\bibitem{Sm1875}
Henry~J.S. Smith.
\newblock On the integration of discontinuous functions.
\newblock {\em Proc. Lond. Math. Soc.}, 6:140--153, 1875.

\bibitem{Can1883}
G.~Cantor.
\newblock \"{U}ber unendliche, lineare punktmannigfaltigkeiten v, [on infinite,
  linear point-manifolds (sets)].
\newblock {\em Math. Ann.}, 21:545--591, 1883.

\bibitem{StSh05}
Elias~M. Stein and Rami Shakarchi.
\newblock {\em Real analysis}.
\newblock Princeton Lectures in Analysis, III. Princeton University Press,
  Princeton, NJ, 2005.

\bibitem{Hos94}
J.~R.~M. Hosking.
\newblock Moments of order statistics of the {C}antor distribution.
\newblock {\em Statist. Probab. Lett.}, 19(2):161--165, 1994.

\bibitem{KnPr96}
Arnold Knopfmacher and Helmut Prodinger.
\newblock Explicit and asymptotic formulae for the expected values of the order
  statistics of the {C}antor distribution.
\newblock {\em Statist. Probab. Lett.}, 27(2):189--194, 1996.

\end{thebibliography}

\end{document}